\documentclass[12pt,thmsa]{article}
\usepackage{amsfonts}
\usepackage{amssymb}
\newtheorem{teo}{Theorem}[section]

\newtheorem{proposition}[teo]{Proposition}

\newtheorem{obs2}[teo]{Remark}

\newtheorem{tea}{Theorem}[subsection]

\newtheorem{no2}[teo]{Note}

\newtheorem{no3}[tea]{Note}

\newcommand{\Gal}{{\rm Gal}}
\newcommand{\Frob}{{\rm Frob }}

\newcommand{\mod}{{\rm mod}}

\newcommand{\Q}{\mathbb{Q}}

\newcommand{\F}{{\mathbb{F}}}
\newcommand{\Z}{{\mathbb{Z}}}

\begin{document}
\title{{\bf On the modularity of rigid Calabi-Yau threefolds: Epilogue
}}
\author{Luis Dieulefait
%\thanks{Supported by .....
%European Research Network ``Arithmetic Geometry"}
\\
Dept. d'\'{A}lgebra i Geometria, Universitat de Barcelona;\\
Gran Via de les Corts Catalanes 585;
08007 - Barcelona; Spain.\\
e-mail: ldieulefait@ub.edu\\
%\thanks{Partially supported     
%by MCYT grant BFM2003-01898}\\   
}
\date{\empty}

\maketitle

\vskip -20mm
%\titlerunning
%

\begin{abstract} In a recent paper of F. Gouvea and N. Yui a detailed account is given of a patching argument due to Serre that proves that the modularity of all rigid Calabi-Yau threefolds defined over $\Q$ follows from Serre's modularity conjecture (now a theorem). In this note we give an alternative proof of this implication. The main difference with Serre's argument is that instead of using as main input residual modularity in infinitely many characteristics we just require residual modularity in a suitable characteristic. This is combined with effective Cebotarev.\\
\end{abstract}

\section{From residual modularity to modularity in rational compatible systems}

 A well-known patching method due to Serre allows to deduce modularity of all rigid Calabi-Yau threefolds over $\Q$ from Serre's modularity conjecture (cf. [Se2] both for the statement of Serre's conjecture and for the patching method of Serre), and in a recent paper by F. Gouvea and N. Yui a detailed account of this method was given (cf. [GY]). In this note we will give an alternative proof of the modularity of these varieties. The main difference with Serre's method is that we will only require residual modularity (i.e., Serre's conjecture) in one suitable characteristic.\\
 
 Recall that in previous works of the author and Manoharmayum, modularity of rigid Calabi-Yau threefolds over $\Q$ was proved under mild local conditions, for example, assuming that the variety has good reduction at $3$ (cf. [DM], [HKS]).\\
 
Before starting, we make an observation that is not necessary but will help to simplify the proof (this observation was also incorporated in a second version of [GY] to simplify their proof): Given the compatible family $\{  \rho_\ell \}$ of Galois representations attached to a rigid Calabi-Yau threefold $X$ over $\Q$, we know that the family is potentially modular (results of R. Taylor, cf. [T]), and as noticed by Taylor a corollary of this is that the family is strictly compatible, meaning that locally at any prime $q$ of bad reduction the Weil-Deligne representation corresponding to the restriction of $\rho_\ell$ to the decomposition group at $q$ is independent of $\ell$. In particular, the conductor of the family is well-defined (for the application we have in mind, uniformly bounded would have been enough). Notice that ``potential modularity" is one of the main tools used to prove Serre's conjecture, so  it seems natural to apply it. This way we will avoid the subtle point of having to select carefully the residual characteristics so that the conductor gets bounded (cf. [Se2]): now we know that the family has a constant conductor.\\
 
Now we want to show modularity of the compatible family applying Serre's conjecture in one characteristic. Our main result is the following:\\

\begin{proposition} Let $X$ be a rigid Calabi-Yau threefold defined over $\Q$. Given the conductor $C$ of the compatible family $\{ \rho_\ell\}$ attached to $X$, there is a bound $B$ (depending on $C$) such that: If we consider $\rho_2$ and $\rho_p$ for a prime $p>B$ we have the following:\\
If $\rho_p$ is residually (irreducible and) modular, then $\rho_2$ is modular.
\end{proposition}

As a corollary of this proposition, since from the work of Khare, Wintenberger, the author and Kisin (cf. [KW1], [KW2], [D] and [K]) we know that Serre's conjecture is true, and since it is well-known that for almost every $p$ the residual representation will be irreducible, we conclude that $\rho_2$ is modular, thus by definition $X$ is modular.\\

Remark 1: the proof that for almost every $p$ the residual mod $p$ representation is irreducible is exactly the same as the one given by Ribet (cf. [Ri]) for the Galois representations attached to cuspidal modular forms. The ingredients are: a compatible family of odd two-dimensional Galois representations of constant conductor $C$ such that for every $p \nmid C$ the $p$-adic representation in the family is known to be crystalline of Hodge-Tate weights $\{0, k-1 \}$ independent of $p$ (in our case $k=4$). From this it can easily be deduced as in [Ri] (using Fontaine-Laffaille theory to deduce that for $p>k$, $p \nmid C$, the residual Serre's weight is $k$, i.e., a uniform description also for the ramification at $p$ of the residual representations) that residual reducibility in infinitely many characteristics implies that the representations in characteristic $0$ are themselves reducible for every $\ell$, in our case: $\rho_\ell \cong \epsilon \oplus \epsilon^{-1} \chi^3$, where $\epsilon$ is a finite order character unramified outside primes in $C$ and $\chi$ is the $\ell$-adic cyclotomic character. But this reducibility contradicts the fact that, as follows from the Weil's conjectures proved by Deligne, these geometric Galois representations are pure, namely, for every $p \nmid C \cdot \ell$, the two eigenvalues of $\rho_\ell (\Frob  \; p)$ have the same absolute value.\\
 
Proof of the Proposition: The tools we use are effective Cebotarev combined with the Weil bounds for the traces $a_q$ of the Galois representations attached to $X$ and to any weight $4$ newform.\\
We know by effective Cebotarev that given a conductor $C$ and two $2$-adic Galois representations $\rho$ and $\rho'$ (odd, two-dimensional, irreducible, of conductor dividing $C$) living both in compatible families which are pure and have the same weights (and thus the same Weil bounds apply to the traces of both), there is a bound $B'$ such that if for every unramified prime $q < B'$ the traces of the two representations at $\Frob \; q$ agree, then the two representations are isomorphic (cf. [Se1], section 8.3, see in particular the comments at the end of the section and the additional remark 636.7 on page 715). In fact, one can take $B'$ to be some suitable power of $C$.\\
 
Now the bound $B$ for the proposition has to be chosen in the following way:\\
$B$ is large enough so that:\\
a) it is larger than any prime of bad reduction of $X$ and larger than $k=4$ (this is needed to ensure that the $p$-adic Galois representation that we consider is crystalline and that its residual Serre's weight can be computed via Fontaine-Laffaille).\\
b) Because of a), for any $p > B$ such that the residual representation is irreducible the residual modularity of the $p$-adic representation attached to $X$ implies the congruence with a cusp form $f$ of weight $4$ and level dividing $C$ (because of Serre's conjecture in its strong form, cf. [Se2]). In the finite dimensional space of cuspidal modular forms of weight $4$ and level $C$ there is a bound $d$ for the degrees of the fields of coefficients of all Hecke eigenforms. Then, for such a $d$, it is easy to see using the Weil bounds for the $a_q$ of $f$ and those of $X$ (they satisfy the same bounds because both are Galois representations pure of Hodge-Tate weights $\{0,  3\}$) that: given the bound $B'$ introduced above (the one for effective Cebotarev), there exists $B$ such that if for some $p > B$ and $q < B'$, $q \nmid C$, we have a congruence
$$ a_q (f) \equiv a_q(X) \quad \pmod{\mathcal{P}} \qquad (*)$$ 
where $\mathcal{P}$ is a prime in the coefficient field of $f$ dividing $p$, then we can conclude that $a_q(f) = a_q(X)$.\\
This is a standard trick, it follows easily from the Weil bounds, observe that the degree $d$ is bounded (just in terms of $C$), and observe that a priori the coefficients of $f$ may not be in $\Q$, but the Weil bounds also apply to their Galois conjugates because they are the coefficients of other modular forms, namely, of the conjugates $f^\sigma$ of $f$. More precisely, we know that:
$$ | a_q (X) | \leq 2 q \sqrt{q}, \qquad | a_q (f) | \leq 2 q \sqrt{q}, \qquad | a_q (f)^\sigma | = | a_q (f^\sigma) | \leq 2 q \sqrt{q} $$
for every $\sigma \in \Gal(\overline{\Q}/ \Q)$. Since the degree of the coefficient field of $f$ is at most $d$, there will be at most $d$ different Galois conjugates of the algebraic integer $a_q(f)$, thus congruence (*) implies that $p$ divides the rational integer:
$$ \prod_{\sigma} ( a_q(f)^\sigma - a_q(X))$$
where $\sigma$ runs in a set of at most $d$ elements. But:
$$ | \prod_{\sigma}  (a_q(f)^\sigma - a_q(X) |  \leq  (4 q \sqrt{q})^{d}$$
Thus if we take $p > (4 q \sqrt{q})^{d}$ we deduce that:
$$ \prod_{\sigma}  (a_q(f)^\sigma - a_q(X) )= 0$$
which implies that $a_q(f) \in \mathbb{Z}$ and $a_q(f) = a_q(X)$. \\
From what we have just seen if follows that taking $B = (4 B' \sqrt{B'})^{d}$ whenever $p$ is bigger than such $B$, for every $q$ of good reduction such that $q < B'$ we have:
$$ p > B = (4 B' \sqrt{B'})^{d} > (4 q \sqrt{q})^{d} $$
and therefore this bound $B$ solves the problem: for any $p$ bigger than $B$ and any $q$ of good reduction smaller than $B'$ congruence (*) modulo $p$ implies that $a_q(f) = a_q(X)$.\\
\newline 
Now we are ready to finish the proof of the proposition. If $p$ is chosen such that  the $\mod \; p$ representation attached to $X$ is irreducible (and thus modular), and $p > B$ for a bound $B$ as in a) and b) above, then the congruence mod $\mathcal{P}$ with some $f$ of weight $4$ and level dividing $C$ implies that the $2$-adic Galois representations attached to $X$ and $f$ have the same coefficients $a_q$ for any unramified $q$ up to $B'$, but this implies (effective Cebotarev) that the two representations are indeed isomorphic, and this concludes the proof.\\

Remark 2: We could have added in the definition of $B'$ the additional condition that it must be  bigger than the Sturm's bound for the space of cuspforms of level $C$ and weight $4$. This way, when we choose $B$ as in the previous proof and we conclude from congruence (*) that $a_q(f) = a_q (X) \in \Z$ for every prime $q \nmid C$, $q < B'$, this implies (by a well-known property of Sturm's bound) that ALL eigenvalues of $f$ are in $\Z$, i.e, that the coefficient field of $f$ is $\Q$. This way we can assume that when we apply effective Cebotarev we are applying it to compare two Galois representations with coefficients in $\Z$. \\  
 
Remark 3: A nice feature of the proof is that given a prime $p_1$ and a sufficiently large prime $p_2$ we are showing that given two geometric Galois representations $\rho_{p_1}$ and $\rho_{p_2}$ that are strictly compatible, then from the residual modularity of $\rho_{p_2}$ we deduce the modularity of $\rho_{p_1}$, thus (by definition) also of $\rho_{p_2}$.\\
 
%It only requires two primes $p_1$ and $p_2$, and it is some sort of ``modularity lifting theorem" in characteristic $p_2$
% that depends just on the existence of the ``fellow" representation for the other prime $p_1$.\\
 
Remark 4:
A third proof of modularity for all rigid Calabi-Yau threefolds defined over $\Q$ (using Serre's conjecture) can be easily deduced from modularity lifting theorems:
It is more or less automatic from a modularity lifting theorem proved by Taylor and (independently) by Diamond-Flach-Guo (cf. [T] and [DFG]) that if $ p> 2k-1 = 7$ is a prime such that: $X$ has good reduction at $p$ and  the $\mod \; p$ representation is irreducible and modular, then $X$ is modular.\\ 
The only technical point to apply the result of [DFG] is the usual condition that the restriction of the residual representation to the quadratic extension ramified only at $p$ has to be irreducible, but using results of Serre and Ribet on dihedral representations (and the description of the action of inertia at $p$ due to Fontaine-Laffaille theory) we easily see that if $p >7$ this case cannot occur.\\
Thus, taking $p$ sufficiently large we deduce modularity of $X$ as a combination of Serre's conjecture with theorems \`{a} la Wiles. Moreover, if we combine with two results of Skinner-Wiles we can remove the condition on residual irreducibility and show that modularity follows just from the existence of a prime of good reduction $p > k-1 = 3$, $p \neq 2k-3 = 5$ and the truth of Serre's conjecture for this $p$: this is exactly how modularity is deduced in [DM] from Serre's conjecture over $\F_7$ (in fact, one can even take $p=3$ if $X$ has good reduction at it, thanks to another  combination of modularity lifting theorems: see my remark at [HKS]).

\section{Bibliography}
%The statement of Serre's conjecture is given in [Se87] and its proof is given in [Di06], [KW06], [KW07] and [Ki06]. The preprint of Gouvea and Yui %mentioned in this note is [GY09]. It is worth recalling that a proof of modularity for rigid Calabi-Yau threefolds defined over $\Q$ under mild %conditions was given in [DM03]. \\

[DFG] Diamond, F., Flach, M., Guo,  L., {\it The Tamagawa number conjecture of adjoint motives of modular forms}, Ann. Sci. Ec. Norm. Sup. {\bf 37} (2004) 663-727
\newline
[Di] Dieulefait, L. V., {\it Remarks on Serre's modularity conjecture}, preprint (2006)
\newline
[DM] Dieulefait, L. V., Manoharmayum, J., {\it Modularity of rigid Calabi-Yau threefolds over $\Q$}, in
{\it Calabi-Yau Varieties and Mirror Symmetry}, Fields Institute Communications Series, AMS {\bf 38} (2003) 159-166
\newline
[GY] Gouvea, F., Yui, N., {\it Rigid Calabi-Yau threefolds over $\Q$ are modular: a Footnote to Serre}, preprint (2009)
\newline
[HKS] Hulek, K., Kloosterman, R., Schuett, M., {\it Modularity of Calabi-Yau varieties}, in
{\it Global Aspects of Complex Geometry}, Springer-Verlag 2006 (F. Catanese, H. Esnault, A. Huckleberry, K. Hulek and T. Peternell, eds.), 271-309
\newline
[KW1] Khare, C., Wintenberger, J.-P., {\it Serre's modularity conjecture (1)}, Inventiones Mathematicae {\bf 178} (2009) 485-504
\newline
[KW2] Khare, C., Wintenberger, J.-P., {\it Serre's modularity conjecture (2)}, Inventiones Mathematicae {\bf 178} (2009) 505-586
\newline
[K] Kisin, M., {\it   Modularity of $2$-adic Barsotti-Tate representations}, Inventiones Mathematicae {\bf 178} (2009) 587-634
\newline 
[R] Ribet, K., {\it On $l$-adic representations attached
 to modular
forms II}, Glasgow  Math.  J.  {\bf 27} (1985) 185-194
\newline
[Se1] Serre, J.-P., {\it Quelques applications du th{\'e}or{\`e}me de densit{\'e} de Chebotarev}, in {\it Oeuvres, Volume III}, Springer-Verlag 1986, 563-641
\newline
[Se2] Serre, J.-P., {\it Sur les repr{\'e}sentations modulaires de degr{\'e}
$2$ de $\Gal(\bar{\mathbb{Q}} / \mathbb{Q})$}, Duke Math. J. {\bf 54}
(1987) 179-230
\newline
[T] Taylor, R., {\it On the meromorphic continuation of degree two
 L-functions}, Documenta Mathematica, Extra Volume: John Coates' Sixtieth Birthday (2006) 729-779

\end{document}